\documentclass{amsproc}
\usepackage{graphicx}
\usepackage{amssymb}
\usepackage{epstopdf}
\usepackage{amsmath,amsthm,amscd,amssymb}
\usepackage{latexsym}
\usepackage[colorlinks,citecolor=red,pagebackref,hypertexnames=false]{hyperref}
\usepackage{geometry}                
\numberwithin{equation}{section}

\theoremstyle{plain}
\newtheorem{theorem}{Theorem}[section]
\newtheorem{lemma}[theorem]{Lemma}

\theoremstyle{definition}
\newtheorem{definition}[theorem]{Definition}

\newtheorem{case[theorem]}{Case}

\theoremstyle{remark}
\newtheorem{remark}[theorem]{Remark}

\def\agtxt{}
\def\od{\mathbb O(d)}
\def\odk{\mathbb O(d-k)}
\def\R{\mathbb R}

\def\bx{{\bf x}}
\def\by{{\bf y}}
\def\stabx{\text{Stab}(\bx)}
\def\staby{\text{Stab}(\by)}

\numberwithin{equation}{section}

\def\bt{\mathbf t}

\begin{document}

\title{A group-theoretic viewpoint on Erd\H os-Falconer problems \\ and the Mattila integral 
}


\author{A. Greenleaf, A. Iosevich, B. Liu and E. Palsson}

\date{today}

\email{allan@math.rochester.edu} 
\email{iosevich@math.rochester.edu}
\email{bliu19@z.rochester.edu}
\email{palsson@math.rochester.edu}

\address{Department of Mathematics, University of Rochester, Rochester, NY 14627}

\thanks{The first two authors were partially supported by  NSF Grants DMS-0853892 and DMS-1045404,  resp.}

\begin{abstract} 
We obtain nontrivial exponents for  Erd\H os-Falconer type point configuration problems. 
Let $T_k(E)$ denote the set of distinct congruent $k$-dimensional simplices determined by $(k+1)$-tuples of points from $E$. 
For $1\le k\le d$, we prove that there exists a $t_{k,d}<d$ such that, if $E \subset {\Bbb R}^d,\, d \ge 2$, with $dim_{{\mathcal H}}(E)>t_{k,d}$, then the ${k+1 \choose 2}$-dimensional Lebesgue measure of $T_k(E)$ is positive.   
Results of this type were previously obtained for triangles in the plane $(k=d=2)$ in \cite{GI12} and 
 for higher $k$ and $d$ in \cite{GGIP12}. 
 We improve upon those exponents,  using a group action perspective, which  also sheds  light on the classical approach to the Falconer distance problem.

\end{abstract} 

\maketitle


\section{Introduction}

\vskip.125in 

One of the most important and far reaching problems in modern geometric measure theory is the Falconer distance problem, which asks: How large does the Hausdorff dimension $s$ of a compact set $E \subset {\Bbb R}^d$, $d \ge 2$, need to be to ensure that the \emph{distance set} of $E$, $\Delta(E):=\{|x-y|: x,y \in E \}\subset \R$, has positive Lebesgue measure? Falconer proved that $s>\frac{d}{2}$ is necessary, up to the endpoint, and conjectured that it  is also sufficient \cite{Falc86}. The best  exponent known to date is $\frac{d}{2}+\frac{1}{3}$, due to Wolff in the plane \cite{W99} and Erdo\~{g}an in higher dimensions \cite{Erd05}. 

A natural extension of the Falconer  distance problem is the \emph{congruent simplex} problem  \cite{EHI12,GI12,GGIP12}.
We say that $\{x^1,\dots, x^{k+1}\}\subset\mathbb R^d$  is \emph{nondegenerate} (or \emph{affinely independent}) if \linebreak$\{x^2-x^1,x^3- x^1,\dots,x^{k+1}-x^1\}$ is linearly independent. This condition  is of course invariant under permutations, and is equivalent with the convex hull of $\{x^1,\dots,x^{k+1}\}$  having positive $k$-dimensional volume and thus being reasonably called the \emph{$k$-simplex generated by} $x^1,\dots,x^{k+1}$, denoted $\Delta(x^1,\dots,x^{k+1})$. 

\vskip.125in

Given a set $E \subset {\Bbb R}^d$, let $E^{k+1}:=E \times E \times \dots \times E $, $(k+1)$-times.

\begin{definition}\label{congruent} 
Let $d\ge2$ and $1\le k\le d$. 
Given a  set $E \subset {\Bbb R}^d$,  define \emph{the set of distinct congruent simplices determined by $E$} to be  $T_k(E):=E^{k+1}/ \sim$, where $(x^1, \dots, x^{k+1}) \sim (y^1, \dots, y^{k+1})$ iff ${\{x^i\}}_{j=1}^{k+1}$ and ${\{y^i\}}_{j=1}^{k+1}$ are nondegenerate and $|x^i-x^j|=|y^i-y^j|$ for $1 \leq i<j \leq k+1$. 
\end{definition} 

{\agtxt There is thus a map $T_k(E)\hookrightarrow{\Bbb R}^{k+1 \choose 2}$, well-defined modulo permutations (which have no effect on positivity of Lebesgue measure and hence will be ignored),
$$\big[(x^1,\dots,x^{k+1})\big]\longrightarrow \Big(|x^i-x^j|\Big)_{1\le i<j\le k+1}.$$}


One  may also consider \emph{similar} simplices instead of congruent ones:
\vskip.125in

\begin{definition}\label{similar} 
Given a compact set $E \subset {\Bbb R}^d$ define  \emph{the set of distinct  similar simplices determined by $E$} to be  $S_k(E):= E^{k+1} / \sim$, where $(x^1, \dots, x^{k+1}) \sim (y^1, \dots, y^{k+1})$ iff ${\{x^i\}}_{j=1}^{k+1}$ and ${\{y^i\}}_{j=1}^{k+1}$ are nondegenerate and, for some $\lambda>0$, $|x^i-x^j|=\lambda |y^i-y^j|$,  for $1 \leq i<j \leq k+1$. \end{definition} 

{\agtxt By  considerations similar to those for $T_k(E)$, one  can   view} $S_k(E)$ as a subset of the projective space ${\Bbb R\Bbb P}^{{k+1 \choose 2}-1}$ or, in local coordinates, 
${\Bbb R}^{{k+1 \choose 2}-1}$.
\vskip.125in

In this paper, we obtain {\agtxt improved (i.e., reduced)  lower bounds on the Hausdorff \linebreak dimension of $E$  that guarantee that $T_k(E)$ and   $S_k(E)$ are  of positive $k+1 \choose 2$ and ${k+1 \choose 2}-1$ \linebreak dimensional Lebesgue measure, resp. The central idea is} a geometric mechanism for studying such problems based on group actions, a method that sheds some new light even on the classical approach to the Falconer distance problem. Our first two results contain the essential features of the method. 

\begin{theorem} \label{method} Let $E$ be a compact set in $ {\Bbb R}^d, \, d \ge 2$, and $\mu$ a finite, nonnegative measure supported on $E$. For  $g\in\od $, the orthogonal group on $\R^d$, define a measure  $\nu_g$, supported  on $E-gE$, by the relation 
\begin{equation} \label{numeasure} \int_{\R^d} f(z)\, d\nu_g(z):=\int_E \int_E f(u-gv)\, d\mu(u)\, d\mu(v),\, f\in C_0(\R^d). \end{equation} 
Define  also a measure $\nu$ on $T_k(E)\subset\R^{k+1\choose 2}$ by 
\begin{equation} \label{kmeasure} \int f(\bt)\, d\nu(\bt)=\int \dots \int f\left(|x^1-x^2|, \dots, |x^i-x^j|,\dots, |x^k-x^{k+1}|\right) d\mu(x^1) \dots d\mu(x^{k+1}), \end{equation} where the entries of the ${k+1 \choose 2}$-vector $\bt$ are the distances $t_{ij}$ from $x^i$ to $x^j$, $1 \leq i<j \leq k+1$.

\vskip.125in 

Then, if $\nu_g$ is absolutely continuous for $a.e. \, g\in\od$, with density also denoted $\nu_g$, and
\begin{equation} \label{knorm} \int_{\od} \int_{\R^d} \nu_g^{k+1}(x)\, dx\,dg<\infty, \end{equation}
where $dg$ is Haar measure on $\od$, 
then the measure $\nu$ in (\ref{kmeasure}) has an $L^2$ density and    ${\mathcal L}^{k+1 \choose 2}(T_k(E))>\nolinebreak0$. 

\end{theorem} 

\vskip.125in 

We  obtain an analogous result for similarity classes.
\vskip.125in

\begin{theorem} \label{methodsimilar} Let $E$ and $\mu$ be as in Thm. \ref{method}. For $a\in\R^+,\, g\in\od$, define 
a measure $\nu_{a,g}$  by
\begin{equation} \label{agmeasure} \int_{\R^d} f(z) \, d\nu_{a,g}(z):=\int_E \int_E f(u-agv)\, d\mu(u)\, d\mu(v), \end{equation} 
Let $I\subset  {\Bbb R}^{+}$ be a a compact interval. Then, if $\nu_{a,g}$ is absolutely continuous for $\hbox{a.e.} \, (a,g)\in I\times\od$ and

\begin{equation} \label{knormsimilar} \int_I \int_{\od} \int_{\R^d} \nu_{a,g}^{k+1}(x)\, dx\, dg\, \frac{da}{a}<\infty, \end{equation} 
then ${\mathcal L}^{{k+1 \choose 2}-1}(S_k(E))>0$. 

\end{theorem} 

\vskip.125in 

As applications of  Theorems \ref{method} and \ref{methodsimilar}, one obtains:

\begin{theorem} \label{mainsimplex} Let $E \subset {\Bbb R}^d$, $d \ge 2$, and $2 \leq k \leq d$. Suppose that 
\begin{equation} \label{exponentcongruent} dim_{{\mathcal H}}(E)>t_{k,d}:= \frac{dk+1}{k+1}.
 \end{equation}
Then ${\mathcal L}^{k+1 \choose 2}(T_k(E))>0$. 
If $d=k=2$, the same conclusion holds if $dim_{{\mathcal H}}(E)>\frac{8}{5}$. 

\vskip.125in 

Now suppose that 

\begin{equation} \label{exponentsimilar} dim_{{\mathcal H}}(E)>s_{k,d} := \frac{dk}{k+1}. 
\end{equation} 
Then ${\mathcal L}^{{k+1 \choose 2}-1}(S_k(E))>0$.  

\end{theorem} 

\vskip.125in 

\begin{remark} In order to illustrate the extent to which the exponents in Theorem \ref{mainsimplex} improve on those in \cite{GI12,GGIP12}, consider the case $k=d$ where $ s_{d,d}=d-1+\frac{2}{d+1}$.
The exponent obtained in \cite{GI12,GGIP12} is $s'_{d,d}=d-\frac{1}{2}+\frac{1}{2d}$. Note that  $s_{d,d}<s'_{d,d}$ for every $d \ge 2$, and, asymptotically  as $d\to\infty$,  the improvement is from $d-\frac{1}{2}$ to $d-1$. To put this in perspective, one has the following lower bound, which shows that for $k=d$ one cannot do better than $d-1$; for further discussion, see Sec. \ref{sec sharpness}.
\end{remark} 

\begin{remark}
While the results of Theorem \ref{mainsimplex} significantly improve and extend the exponents in \cite{EHI12,GI12,GGIP12},  the group-theoretic nature of our methods also casts new light upon the classical Mattila integral (see Sec. \ref{sec free}), potentially leading to further progress on related problems. 
\end{remark}

\vskip.125in 

\begin{theorem} \label{sharpness} Let $\alpha_{k,d}$ denote the optimal exponent for the congruent $d$-dimensional simplex problem, i.e., $\alpha_{k,d}$ is the infimum of those $\alpha$ for which ${\mathcal L}^{k+1 \choose 2}(T_k(E))>0$ whenever $dim_{{\mathcal H}}(E)>\alpha$. Then 
$$\alpha_{k,d} \ge \max \left\{k-1, \frac{d}{2} \right\}.$$
Moreover,  $\alpha_{2,2} \ge \frac{3}{2}$. 
\end{theorem} 

\vskip.125in

We thank Ciprian Demeter for pointing out the relevance of \cite{BD13} to the question of  sharpness examples, and an anonymous referee for suggesting many improvements to the paper.

\vskip.25in

\section{Proofs of Theorem \ref{method} and Theorem \ref{methodsimilar}} 
\label{geometricapproach}

Motivated by the geometric viewpoint in \cite{IL13}, the essence of our approach is the following. Define a measure $d\nu$ on $\R^{k+1\choose 2}=\R^{\frac{k(k+1)}{2}}$, with support in $T_k(E)$, as in (\ref{kmeasure}) above.   We will show that  to prove Thm. \ref{method} it  suffices to obtain an upper bound on the $L^2$ norm of the density, i.e., the Radon-Nikodym derivative of  $d\nu$, which we denote by $\nu(\bt)$. 
We start by showing that

\begin{eqnarray}\label{thickened}
\int \nu^2(\bt)\, d\bt\le c_{k,d}\cdot
 \liminf_{\epsilon \to 0}& \epsilon^{- \frac{k(k+1)}{2}} \mu^{2(k+1)} \big\{ &(x^1, \dots, x^{k+1}, y^1, \dots, y^{k+1})\in(\R^d)^{2(k+1)}:\nonumber \\
 &  &\left||x^i-x^j|-|y^i-y^j|\right|\leq \epsilon, \ 1 \leq i<j \leq k+1 \big\},
\end{eqnarray} 
where $\mu^{2(k+1)}$ denotes $\mu \times \dots \times \mu$, $2(k+1)$ times, with the proof showing that if the RHS of (\ref{thickened}) is finite,
then in fact $d\nu$ is absolutely continuous with respect to Lebesgue measure $d\bt$, with density $\nu(\bt)\in L^2$ . 
Let $\phi\in C_0^\infty(\R^{\frac{k(k+1)}2}),\, \phi\ge 0,\, supp(\phi)\subset \{\bt|\le 1\},\, \int \phi \, d\bt=1$, and $\phi_\epsilon(\cdot)=\epsilon^{-\frac{k(k+1)}2}\phi(\epsilon^{-1}\cdot),\, 0<\epsilon<\infty$, the resulting approximate identity. 
Setting $\nu_\epsilon=\phi_\epsilon * d\nu\in C_0^\infty$, one has $d\nu=wk^*\!-\!\lim_{\epsilon\to 0} \nu_\epsilon$, and
(\ref{thickened}) will follow if one shows that $\liminf_{\epsilon\to 0} ||\nu_\epsilon||_{L^2}^2=C<\infty$.

\vskip.125in 

Now, 
$$\nu_\epsilon(\bt)=\langle d\nu,\phi_\epsilon(\cdot-\bt)\rangle = \int \phi_\epsilon\left(\left( |x^i-x^j|-t_{ij}\right)_{1\le i<j\le k+1}\right)\, d\mu(x^1)\cdots d\mu(x^{k+1}).$$
Due to the nonnegativity of $\phi_\epsilon$ and $d\mu$, this is dominated by
$$\int \prod_{1\le i<j\le k+1}\, \epsilon^{-1}\chi\left\{\left||x^i-x^j|-t_{ij}\right|<\epsilon\right\}\, d\mu(x^1)\cdots d\mu(x^{k+1}),$$
where $\chi A(\cdot)$ denotes the characteristic function of a set $A$, and thus

\begin{eqnarray}\label{squared}\nonumber
||\nu_\epsilon||_{L^2}^2\lesssim \int &  & \prod_{1\le i<j\le k+1} 
\epsilon^{-1}\chi\left\{\left||x^i-x^j|-t_{ij}\right|<\epsilon\right\}
\prod_{1\le i<j\le k+1} 
\epsilon^{-1}\chi\left\{\left||y^i-y^j|-t_{ij}\right|<\epsilon\right\} \\
& & \qquad d\mu(x^1)\cdots d\mu(x^{k+1})\, d\mu(y^1)\cdots d\mu(y^{k+1})\,d\bt.
\end{eqnarray}
Now, by the triangle inequality, one has
$$\chi\left\{\left||x^i-x^j|-t_{ij}\right|<\epsilon\right\}\cdot\chi\left\{\left||y^i-y^j|-t_{ij}\right|<\epsilon\right\}\le
\chi\left\{\left||x^i-x^j|-|y^i-y^j|\right|<2\epsilon\right\},$$
and thus, integrating out $d\bt$, the RHS of (\ref{squared}) is
$$\lesssim \epsilon^{-\frac{k(k+1)}2}\int \prod_{1\le i<j\le k+1} 
\chi\left\{\left||x^i-x^j|-|y^i-y^j|\right|<2\epsilon\right\} 
d\mu(x^1)\cdots d\mu(x^{k+1})\, d\mu(y^1)\cdots d\mu(y^{k+1}).
$$
Taking the $\liminf$ as $\epsilon\to 0$ yields the RHS of (\ref{thickened}).

\vskip.125in

To continue, we next introduce some notation. We denote an ordered $(k+1)$-tuple of elements of $\R^d$ by $\bx:=(x^1,\dots,x^{k+1})$.
If the corresponding set $\{x^1,\dots,x^{k+1}\}$ is nondegenerate (i.e., affinely independent), then
$$\pi(\bx):=span\{x^2-x^1,\dots, x^{k+1}-x^1\}$$
is a $k$-dimensional linear subspace of $\R^d$. Let $\Delta(\bx)$ be the (unoriented) simplex  generated by $\{x^1,\dots,x^{k+1}\}$, i.e., the closed convex hull,  which is contained in the affine plane $x^1+\pi(\bx)$. Both $\pi(\bx)$ and $\Delta(\bx)$ are independent of the order of the $x^j$. 
If $\{y^1,\dots, y^{k+1}\}$ is congruent to $ \{x^1,\dots,x^{k+1}\}$, as defined in Def. \ref{congruent}, then an elementary argument shows that, up to permutation of $y^1,\dots, y^{k+1}$, there exists a $g\in\od $ such that $x^j-x^1=g(y^j-y^1),\, 2\le j\le k+1$, which is equivalent with $x^j-x^i=g(y^j-y^i), 1\le i<j\le k+1$, and $\Delta(\bx)=(x^1-gy^1)+ g\Delta(\by)$ . The group $\od$ acts on the Grassmanians $G(k,d)$ and $G(d-k,d)$ of $k$ (resp., $d-k$) dimensional linear subspaces of $\R^d$, and if $\bx$ is congruent to $\by$, one has
$\pi(\bx)=g\pi(\by)$ and $\pi(\bx)^\perp=g\big(\pi(\by)^\perp\big)$.
The set of $g\in\od$ fixing $\Delta(\bx)$ is a conjugate of $\odk\subset\od$, and we refer to this as the \emph{stabilizer} of $\Delta(\bx)$, denoted $\stabx$.

\vskip.125in

For $\bx,\, \by$ congruent as above, let $\tilde{g}\in\od$  be such that $\pi(\bx)=\tilde{g}\pi(\by)$.
Then,  $x^i-x^j=\tilde{g}h(y^i-y^j)$ for all $h\in \staby$. For each $\by$, take a  cover of $\od / \staby$ by balls  of radius $\epsilon$ (with respect to some Riemannian metric) with finite overlap.  
Since the dimension of
 $\od/\staby$ is that of  $\od / \odk$, namely
$$\frac{d(d-1)}{2} - \frac{(d-k)(d-k-1)}{2} = kd - \frac{k(k+1)}{2},$$
we need $N(\epsilon)\sim C \epsilon^{-\left(kd - \frac{k(k+1)}{2}\right)}$ balls to cover it.
Choose  sample points, $\tilde{g}_m(\by),\, 1\le m\le N(\epsilon)$, one in each of the balls.

From basic geometry one sees that the set
$$ \left\{(\bx, \by): \left||x^i-x^j|-|y^i-y^j|\right|\leq \epsilon, \ 1 \leq i<j \leq k+1 \right\} $$
is contained in
$$ \bigcup\limits_{m=1}^{N(\epsilon)} \left\{(\bx, \by): \left|(x^i-x^j)-\tilde{g}_m(\by) h(y^i-y^j)\right|\leq C \epsilon, \,\forall\, 1 \leq i<j \leq k+1, h\in\staby \right\} , $$
where $C=2\max\left\{\text{diam}(E),1 \right\}$. 
Thus,  the expression within the $\liminf$ on the RHS of (\ref{thickened}) is bounded above by
$$
\epsilon^{- \frac{k(k+1)}{2}} \sum\limits_{m=1}^{N(\epsilon)} \mu^{2(k+1)} \left\{(\bx,\by): \left|(x^i-x^j)-\tilde{g}_m(\by) h(y^i-y^j)\right|\leq C \epsilon, \, \forall\, 1\, \leq i<j \leq k+1, h\in\staby \right\}
$$
which can also be written as
\begin{multline} \label{gettingthere}
\epsilon^{- kd} \sum\limits_{m=1}^{N(\epsilon)} \epsilon^{kd- \frac{k(k+1)}{2}}\mu^{2(k+1)} \{(\bx,\by): \left|(x^i-\tilde{g}_m(\by)  h y^i)-(x^j-\tilde{g}_m(\by) h y^j)\right|\leq C \epsilon, \\ \forall\, 1 \leq i<j \leq k+1, h\in\staby \} .
\end{multline}
Since this holds for any choice of  sample points $\tilde{g}_m(\by)$, we can  pick these points such that they minimize (up to a factor of 1/2, say) the quantity
$$ \mu^{2(k+1)} \{(\bx,\by): \left|(x^i-\tilde{g}_m(\by) h y^i)-(x^j-\tilde{g}_m(\by) h y^j)\right|\leq \epsilon, \,\forall\, 1 \leq i<j \leq k+1, h\in\staby \} .$$

Now consider the $N(\epsilon)$ preimages, under the natural projection from $\od$, of the balls  used to cover $\od / \staby$; we can label  these $\epsilon$-tubular neighborhoods of the preimages of the sample points  $\tilde{g}_m(\by)$  as $T_1^{\epsilon},\ldots,T_{N(\epsilon)}^{\epsilon}$.  
Since $dim(\od / \staby)=kd - \frac{k(k+1)}{2}$, each  $T^\epsilon_m$ has volume  \linebreak$\sim \epsilon^{kd-\frac{k(k+1)}2}$. The inf over a set is less than or equal to the average over the set, so we obtain
\begin{multline*}
\mu^{2(k+1)} \{(\bx,\by): \left|(x^i-\tilde{g}_m(\by) h y^i)-(x^j-\tilde{g}_m(\by) h y^j)\right|\leq \epsilon, \,\forall\, 1 \leq i<j \leq k+1, h\in\staby \} \\
\lesssim \frac{1}{\epsilon^{kd - \frac{k(k+1)}{2}}} \int\limits_{T_m^{\epsilon}} \mu^{2(k+1)} \{(\bx,\by): \left|(x^i-g y^i)-(x^j-g y^j)\right|\leq \epsilon, \ 1 \leq i<j \leq k+1 \} \, dg.
\end{multline*}
We can thus bound (\ref{gettingthere}) above by
$$\epsilon^{-kd} \sum\limits_{m=1}^{N(\epsilon)} \, \int\limits_{T_m^{\epsilon}} \mu^{2(k+1)} \{(\bx,\by): \left|(x^i-g y^i)-(x^j-g y^j)\right|\leq \epsilon, \ 1 \leq i<j \leq k+1 \} dg. $$
Since the cover has finite overlap, this in turn can be bounded above, up to a  constant $c_{k,d}$,  by
$$ \epsilon^{-kd} \int \mu^{2(k+1)} \{(\bx,\by): \left|(x^i-g y^i)-(x^j-g y^j)\right|\leq \epsilon, \ 1 \leq i<j \leq k+1 \} dg,$$ and taking the liming, we obtain a constant multiple of the expression (\ref{knorm}). This completes the proof of Theorem \ref{method}. 

\vskip.125in

If one  studies \emph{similar} simplices instead of congruent ones, letting $S_k(E)$ as in  Def. \ref{similar}, then the preceding analysis goes through essentially unchanged, except that in place of (\ref{knorm}) we have (\ref{knormsimilar}). This establishes Theorem \ref{methodsimilar}. 

\vskip.25in

\section{Proof of Theorem \ref{mainsimplex}} \label{sec mainproof}

\vskip.125in 

The matters have been reduced in the introduction to the estimation of (\ref{knorm}). We shall need the following result. 

\begin{theorem} \label{wolfferdogan} Let $\mu$ be a compactly supported Borel measure. Then, for $s\ge \frac{d}2$, $\epsilon>0$,
$$ \int_{S^{d-1}} {|\widehat{\mu}(t \omega)|}^2 d\omega \leq C_\epsilon I_s(\mu)\,  t^{\epsilon-\gamma_s},$$  with $\gamma_s=\frac{d+2s-2}{4}$ if $ \frac{d}{2} \leq s \leq \frac{d+2}{2}$, and $\gamma_s=s-1$ for $ s \ge \frac{d+2}{2}$. 
\end{theorem} 

For $s\le\frac{d+2}2$, this is due to Wolff \cite{W99}($d=2$) and Erdo\~{g}an \cite{Erd05} ($d\ge 3$); the easier case of 
$s\ge\frac{d+2}2$ is due to Sj\"{o}lin \cite{Sj93}.

As we note above, the proof of Theorem \ref{mainsimplex} is reduced to the verification of (\ref{knorm}). Let $\psi$ be a smooth cutoff function supported in $\left\{\xi \in {\Bbb R}^d: \frac{1}{2} \leq |\xi| \leq 4 \right\}$ and identically equal to $1$ in $\left\{\xi \in {\Bbb R}^d: 1 \leq |\xi| \leq 2  \right\}$. Let $\nu_{g,j}$ denote the $j$th Littlewood-Paley piece of $\nu_g$, defined by the relation $\widehat{\nu}_{g,j}(\xi)=\widehat{\nu}_g(\xi) \psi(2^{-j}\xi)$. Since $\nu_g$ is compactly supported, we may assume that $j\ge 0$. Using the Littlewood-Paley decomposition of $\nu_{g}^{}$, the integral in (\ref{knorm}) equals
$$ \int \sum_{j_1, \dots, j_{k+1}} \nu_{g,j_1}(x) \nu_{g,j_2}(x) \cdots \nu_{g,j_{k+1}}(x)\,  dx .$$

We can split this sum up into $k(k+1)$ sums of the type where we sum up over indices where $j_1\geq j_2 \geq j_3,\ldots,j_{k+1}$ and permutations thereof. It suffices to show bounds for one of those sums so without loss of generality we may assume we are in the case $j_1\geq j_2 \geq j_3,\ldots,j_{k+1}$. Passing to the Fourier side we can write
$$\int \sum_{j_1 \geq j_2 \geq j_3 \dots, j_{k+1}} \nu_{g,j_1}(x) \nu_{g,j_2}(x) \cdots \nu_{g,j_{k+1}}(x)\, dx$$ 
$$ =\sum_{j_1 \geq j_2 \geq j_3 \dots, j_{k+1}} \int \widehat{\nu}_{g,j_1}*\widehat{\nu}_{g,j_3}*\dots*\widehat{\nu}_{g,j_{k+1}}(\xi) \cdot \widehat{\nu}_{g,j_{2}}(\xi) \, d\xi.$$ 
Now $\widehat{\nu}_{g,j_1}*\widehat{\nu}_{g,j_3}*\dots*\widehat{\nu}_{g,j_{k+1}}$ is supported on scale $2^{j_1}+2^{j_3}+\ldots+2^{j_{k+1}}\sim 2^{j_1}$ while $\widehat{\nu}_{g,j_{2}}$ is supported on scale $2^{j_2}$ so by Plancherel it is clear that the sum vanishes if $j_1-j_{2}>2$. Thus it suffices to consider the case $j_1 = j_2$ and to study
\begin{equation}\label{twoequal}
\sum_{j_1 = j_2 \geq j_3 \dots, j_{k+1}} \int  \nu_{g,j_1}(x) \nu_{g,j_2}(x) \cdots \nu_{g,j_{k+1}}(x) \, dx .
\end{equation}

As  above,  $ \nu_{g,j}(x)=\mu_j*\mu_j(g \cdot)$, so 
$$ {||\nu_{g,j}||}_{\infty} \leq {||\mu_j||}_1 \cdot {||\mu_j||}_{\infty} \leq C2^{j(d-s)},$$ for any $s<dim_{{\mathcal H}}(E)$, since $\mu$ is a Frostman measure supported on $E$ (see e.g. \cite[Chap. 8]{M95}). To see this, observe that ${||\mu_j||}_1 \leq 1$ trivially since $\mu$ is a probability measure and 
$$ |\mu_j(x)|=2^{dj} \left| \mu*\widehat{\psi}(2^j \cdot)(x) \right|$$
$$ \leq C_N 2^{dj} \int {(1+2^j|x-y|)}^{-N} d\mu(y) \leq C'_N2^{j(d-s)}$$ since $\mu$ is a Frostman measure on $E$. Using this estimate on the terms corresponding to the indices $j_3,\ldots,j_{k+1}$ we can bound (\ref{twoequal}) above, up to a fixed constant, by
$$ \sum\limits_{j}\left(\sum\limits_{j_3,\ldots,j_{k+1}\leq j} 2^{(j_3+\ldots+j_{k+1})(d-s)} \right)\nu^2_{g,j}(x) $$
$$ \lesssim \sum\limits_{j} 2^{j(k-1)(d-s)} \nu^2_{g,j}(x) .$$
It follows that we can bound (\ref{knorm}) by a finite sum of terms of the type
\begin{equation} \label{keepitcoming} \sum_j 2^{j(k-1)(d-s)} \cdot \int \int \nu^2_{g,j}(x)\, dx\, dg. \end{equation}

By Plancherel (see the discussion in Sec. \ref{sec free} below), 
$$ \int \int \nu^2_{g,j}(x) \, dx \, dg \approx \int_{2^j}^{2^{j+1}} {\left( \int_{S^{d-1}} {|\widehat{\mu}(t \omega)|}^2 d\omega \right)}^2 t^{d-1} dt$$ 
$$ \leq C'' 2^{j(d-s)} 2^{-j \gamma(s,d)},$$ with the  inequality following from Theorem \ref{wolfferdogan}. 

Let us first handle the case $d \ge 3$. Inserting the last inequality back into (\ref{keepitcoming}), we see that geometric series converges if $(d-s)k-(s-1)<0$, which yields the condition $s>\frac{dk+1}{k+1}$, as claimed. If $d=k=2$, $\gamma(s,2)=\frac{s}{2}$ and the geometric series converges if $s>\frac{8}{5}$. This completes the proof of the first part of Theorem \ref{mainsimplex}. 

To prove the second part, as  explained in Sec. \ref{geometricapproach}, it  suffices to estimate (\ref{knormsimilar}). Following the proof of the first part of Theorem \ref{mainsimplex} above, the second part would follow from the estimate 
\begin{equation} \label{lafa} \left| \int_1^2 \int {|\widehat{\mu}(ag\xi)|}^2 \, dg\, \frac{da}{a}  \right| \leq C{|\xi|}^{-s}, \end{equation} where the reduction to $a \in [1,2]$ is accomplished by simple pigeon-holing and scaling. Indeed, recall that $\mu$ is a Frostman measure supported on $E$ (see e.g. \cite{W03}, Chapter 8) means that for any $\epsilon>0$ there exists $C_{\epsilon}>0$ such that if $B_{\delta}$ is a ball of radius $\delta$ centered at the origin, then 
\begin{equation} \label{frostmanbound} \mu(B_{\delta}) \leq C_{\epsilon} \delta^{s-\epsilon}, \end{equation} where $s$ is the Hausdorff dimension of $E$ and $\epsilon>0$ is arbitrarily small. If $\delta$ is sufficiently small, this quantity is $<\frac{1}{2}$, so the intersection $E$ and the complement of $B_{\delta}$ has $\mu$-measure $>\frac{1}{2}$. Renaming this intersection as $E$ and rescaling, the procedure that does not affect whether the Lebesgue measure of $T_k(E)$ (or $S_k(E)$) is positive, we achieve the desired setup. 

By the action of the orthogonal group on the sphere, (\ref{lafa}) would follow from 

$$ \left| \int_1^2 \int_{S^{d-1}} {\left|\widehat{\mu}\left(a|\xi| \omega \right)\right|}^2\,  d\omega\, \frac{da}{a}  \right| \leq C{|\xi|}^{-s}$$ for any $s<dim_{{\mathcal H}}(E)$. This in turn is proven by observing that 

\begin{eqnarray*}
 \int_1^2 \int_{S^{d-1}} {|\widehat{\mu}(a|\xi| \omega)|}^2\, d\omega \, \frac{da}{a} &=&\int_{|\xi|}^{2 |\xi|} \int_{S^{d-1}} {|\widehat{\mu}(a \omega)|}^2\, d\omega\,  \frac{da}{a} \\
& \leq&  {|\xi|}^{-d}  \int_{|\xi|}^{2 |\xi|} \int_{S^{d-1}} {|\widehat{\mu}(a \omega)|}^2 \, d\omega\, a^{d-1}\,  da\\
&=&{|\xi|}^{-d} \int_{|\xi| \leq |x| \leq 2|\xi|} {|\widehat{\mu}(x)|}^2\,  dx \leq C{|\xi|}^{-s}.
\end{eqnarray*}

\vskip.25in

\section{Sharpness of lower bounds} \label{sec sharpness}

We now turn to the proof of Thm. \ref{sharpness}.
 If $E$ is contained in a $(k-1)$-dimensional plane, every simplex with $k+1$ points in $E$ is degenerate,
so the restriction $\alpha_{k,d}>k-1$ is clear. The lower bound $\alpha_{k,d}>\frac{d}{2}$ is also necessary, since this threshold is needed  to ensure that $\mathcal L^1(\Delta(E))>0$ for  general $E$, as was noted in the introduction. 

The restriction $\alpha_{2,2}>\frac{3}{2}$ follows from a lattice construction and simple number-theoretic analysis, given below but previously obtained by Burak Erdo\~{g}an and the second listed author \cite{EI11}. 
Start by considering the following construction in $\R^d$ for general $d\ge 1$. 
 Let $q_1=2$ and recursively choose $q_i\in\mathbb N$ with $q_{i+1}>q_i^i, \,\forall i\ge 1$. Fix $s<d$ and let $E_i$ denote the $q_i^{-\frac{d}{s}}$-neighborhood of $ {\Bbb Z}^d \cap {[0,q_i]}^d$, scaled by $q_i^{-1}$ so as to be a subset of $[0,1]^d$. Define $E=\cap_i E_i$. The proof that $dim_{{\mathcal H}}(E)=s$ can be found in \cite[Chap. 8]{Falc86} in the  case of $d=1$. The higher dimensional argument follows from the same argument. 

Now let $d=2$. To show that  $\alpha_{2,2}\ge \frac{3}{2}$, let $q=q_i$, for $i$ very large. We claim that 
 
\begin{equation} \label{trianglekaif}
 {\mathcal L}^3(T_2(E_q)) \lesssim q^{-\frac{6}{s}} \cdot \# T_2(\{ {\Bbb Z}^2 \cap {[0,q]}^2 \}) \leq C_{\epsilon} q^{-\frac{6}{s}} \cdot q^{4+\epsilon}. \end{equation} 
In fact, note that, by translation invariance,  in order to count congruence classes determined by the unrescaled ${\Bbb Z}^2 \cap {[0,q]}^2$ it is enough to place one vertex at the origin. Call the remaining vertices $v$ and $w$ and let $|v|=a$, $|w|=b$. For our purposes it is sufficient to know that the number of choices is $\leq Cq^2$. To see this simply observe that squares of the distances from the origin are integers in $[0,2q^2]$, so there cannot possibly be more than $2q^2$ of them. 
Once $|v|$ and $|w|$ are fixed, it remain to compute how many possibilities there are for $|v-w|$. This number cannot exceed the product of the number of integer points on $\{x: |x|=a \}$ and the number of integer points on $\{x: |x|=b \}$. It is well known  that for any $\epsilon>0$, the number of lattice points on the circle of radius $r$ in the plane does not exceed $C_{\epsilon}r^{\epsilon}$; see, e.g., \cite{Lan69}. The estimate (\ref{trianglekaif}) thus follows and we conclude that ${\mathcal L}^3(T_2(E))$ is not in general positive if the Hausdorff dimension of $E$ is smaller than $\frac{3}{2}$. 

Note that  this argument does not show what happens at $s=\frac{3}{2}$. 

\vskip.125in 

Examination of $\alpha_{3,3}$ leads to an interesting lattice point problem.  Take $d=3$ in the construction above. Once again, we place one vertex at the origin and call the remaining vertices $v^1, v^2, v^3$. There are $\approx q^2$ choices for $|v^i|$. It remains to count the number of non-congruent configurations that $v^i$s can form. Each $v^i$ lies on a sphere of radius at most $q$ and it is well-known that the number of lattice points a sphere of radius $r$ in ${\Bbb R}^3$ is $\lessapprox q$. It follows by trivial counting that the number of non-congruent configurations of $v^j$s is $\lessapprox q^3$. We deduce that 

$$ {\mathcal L}^6(T_4(E_q)) \lessapprox q^{-\frac{18}{s}} \cdot q^6 \cdot q^3, $$ which results in the trivial restriction $s>2$. So the question of whether we can obtain a tighter restriction on $dim_{{\mathcal H}}(E)$, needed to ensure that ${\mathcal L}^6(T_4(E))>0$, comes down to estimating the size of the discrete set 

$$T_3(S_1 \cap {\Bbb Z}^3,S_2 \cap {\Bbb Z}^3,S_3 \cap {\Bbb Z}^3),$$ the number of non-congruent triangles with vertices at lattice points on spheres $S_1,S_2,S_3$ of radii $\approx q$. Any estimate of the form 

\begin{equation} \label{3spheres} \# T_3(S_1 \cap {\Bbb Z}^3,S_2 \cap {\Bbb Z}^3,S_3 \cap {\Bbb Z}^3) \leq Cq^{3-\delta} \end{equation} for some $\delta>0$ would immediately allow one to conclude that 
$$ \alpha_{3,3} \ge 2+\delta'$$ for some $\delta'>0$.  We do not know whether (\ref{3spheres}) holds, and pose this question as an open problem that is interesting in its own right. \footnote{After submission of this paper, this question was  answered in the negative by  Demeter \cite{D13}.}

\vskip.25in 

\section{A stationary phase-free proof of the $\frac{d+1}{2}$ exponent in the Falconer problem} \label{sec free}

The purpose of this section is to make a couple of simple observations regarding the Falconer distance conjecture and the methods of proof that have been employed to attack it. First, we apply the results of Sec. \ref{sec mainproof} to the case $k=1$, corresponding to the  Falconer distance problem. Applying (\ref{numeasure}) with $f(z)=e^{-2 \pi i z \cdot \xi}$, we obtain 
$$ \widehat{\nu}_g(\xi)=\widehat{\mu}(\xi) \widehat{\mu}(g \xi),$$ which means that, via Plancherel, the expression in (\ref{knorm}), with $k=1$ is equal to
$$ \int_{\R^d}{|\widehat{\mu}(\xi)|}^2 \left\{ \int_{\od}{|\widehat{\mu}(g \xi)|}^2 dg \right\} d\xi.$$ 

A moment's reflection shows that this quantity equals a constant multiple of 
\begin{equation} \label{mattilakrutoi} \int {\left( \int_{S^{d-1}} {|\widehat{\mu}(t \omega)|}^2 d\omega \right)}^2 t^{d-1} dt, \end{equation} the classical Mattila integral derived in \cite{Mat87}, which has so far been the main tool in the study of the Falconer distance problem. The fact that the boundedness of this integral implies a lower bound on the Lebesgue measure of the distance set is typically derived using the method of stationary phase (see also \cite[Chap. 9]{W03}), but the argument above shows that a group-theoretic argument can be used instead. 

We now establish the fact that the threshold $\frac{d+1}{2}$ for the Falconer distance conjecture can be established using our geometric methods without the use of the method of stationary phase. See also Mitsis \cite{Mi02} for another geometric argument in the context of the Falconer distance problem. 

The argument culminating in (\ref{mattilakrutoi}) above, reproves the classical result due to Mattila, namely that if $E$ is a compact subset of ${\Bbb R}^d$ of Hausdorff dimension $s>\frac{d}{2}$, and the Mattila integral given by (\ref{mattilakrutoi}) is bounded for some Borel measure $\mu$ supported on $E$, then the Lebesgue measure of $\Delta(E)$ is positive. The fact that the expression in (\ref{mattilakrutoi}) is bounded if the Hausdorff dimension of $E$ is greater than $\frac{d+1}{2}$ for any Frostman measure $\mu$ supported on $E$ (see \cite[p.\! 112]{M95} for background on Frostman's Lemma) follows immediately from the following simple observation. Recall the $s$-energy integral of $\mu$, 

\begin{equation} \label{energyintegral} I_s(\mu)=\int \int {|x-y|}^{-s}\, d\mu(x) \, d\mu(y). \end{equation} 

\begin{lemma} \label{trivialspherical} \cite[p.\! 61]{W03} Let $\mu$ be a compactly supported Borel measure on ${\Bbb R}^d,\,d \ge 2$. Then, for any $s \ge \frac{d}{2}$, 
$$ \int_{S^{d-1}} {|\widehat{\mu}(t \omega)|}^2 d\omega \leq CI_s(\mu) t^{-(s-1)},\quad \, 0<t<\infty. $$ 
\end{lemma} 

We give the proof of Lemma \ref{trivialspherical}  for the sake of completeness. 
Let $\phi$ be a radial smooth function with compact support whose Fourier transform is $\ge 1$ on the support of $\mu$. Then it suffices to estimate 
\begin{eqnarray*}
 \int {|\phi*\hat{\mu}(t\omega)|}^2 d\omega &\leq&  \Big(\int \big(\int |\phi(x-t\omega)|^2 {|\widehat{\mu}(x)|}^2 \, d\omega\big)^\frac12\, dx\Big)^2  \\
&\leq& C't^{-(d-1)} \int_{||x|-t| \leq C''} {|\widehat{\mu}(x)|}^2 dx \\
&\leq& C'''I_s(\mu) t^{-s+1},
\end{eqnarray*}  
finishing the proof.

To see how Lemma \ref{trivialspherical} implies the $\frac{d+1}{2}$ exponent for the Falconer problem, it is enough to prove that the Mattila integral (\ref{mattilakrutoi}) is bounded if $\mu$ is a Borel measure supported on a set of Hausdorff dimension greater than $\frac{d+1}{2}$. Using Lemma \ref{trivialspherical} we see that (with $I_s(\mu)$ as in (\ref{energyintegral}) above). 

$$ \int {\left( \int_{S^{d-1}} {|\widehat{\mu}(t \omega)|}^2 d\omega \right)}^2 t^{d-1}dt$$ 

$$ \leq C \int \int t^{d-1} t^{-s+1} {|\widehat{\mu}(t\omega)|}^2 d\omega \, dt$$ 

$$=C \int {|\widehat{\mu}(\xi)|}^2 {|\xi|}^{-s+1} d\xi$$ 

$$=C \int {|\widehat{\mu}(\xi)|}^2 {|\xi|}^{-d+(d-s+1)} d\xi \leq C'' I_s(\mu)$$ if $s>\frac{d+1}{2}$, as desired.

\bigskip

\end{document}